\documentclass[12pt]{article}
\usepackage{amsmath}
\usepackage{epsfig}
\usepackage{amsfonts}
\usepackage{mathrsfs}
\usepackage{latexsym}
\usepackage{graphicx}
\def\R{\mathbf{R}}
\def\C{\mathbf{C}}
\def\bC{\mathbf{\overline{C}}}
\def\Ima{\mathrm{Im}\ }
\def\Rea{\mathrm{Re}\ }
\def\D{{\mathscr{D}}}
\begin{document}
\title{An extremal problem for a class of entire
functions of exponential type}
\author{Alexandre Eremenko\thanks{Supported by NSF grant DMS-0555279.}
$\;$ and Peter Yuditskii\thanks{Supported by Austrian Fund FWF P20413-N18.}}
\date{May 29, 2008}
\maketitle

\begin{abstract}
Let $f$ be an entire function of exponential type
whose indicator diagram is contained in the segment
$[-i\sigma,i\sigma],\sigma>0.$ Then the upper density
of zeros of $f$ is at most
$c\sigma$ where
$c\approx 1.508879$ is the positive solution
of the equation
$$\log(\sqrt{c^2+1}+c)=\sqrt{1+c^{-2}}.$$
This estimate is exact. 
\end{abstract}

We consider the class $E_\sigma,\;\sigma>0$ of
entire functions of exponential
type whose indicator diagram
is contained in a segment
$[-i\sigma,i\sigma],$ which means that
\begin{equation}
\label{0}
h(\theta):=\limsup_{r\to+\infty}\frac{\log|f(re^{i\theta})|}{r}\leq \sigma|\sin\theta|,
\quad|\theta|\leq\pi.
\end{equation}
An alternative characterization of such functions
follows from a theorem of P\'olya \cite{L}:
$$f(z)=\frac{1}{2\pi}\int_\gamma F(\zeta)e^{-i\zeta z}d\zeta,$$
where $F$ is an analytic function in $\bC\backslash[-\sigma,\sigma],\;
F(\infty)=0$, and $\gamma$ is a closed contour
going once around the
segment $[-\sigma,\sigma]$. In other words, the class of entire functions
satisfying (\ref{0}) consists of Fourier transforms of hyperfunctions
supported by $[-\sigma,\sigma]$, see, for example,
\cite{EN} and
\cite{Hor}.

Let $n(r)$ be the number of zeros of $f$ in the
disc $\{z:|z|\leq r\}$,
counting multiplicity. We are interested in the
{\em upper density}
\begin{equation}
\label{1}
D=\limsup_{r\to\infty}\frac{n(r)}{r}.
\end{equation}

If $f$ satisfies the additional condition
\begin{equation}
\label{2}
\int_{-\infty}^\infty\frac{\log^+|f(x)|}{1+x^2}dx<\infty,
\end{equation}
then the limit (density) in (\ref{1}) exists and equals
$\displaystyle (2\pi)^{-1}\int_{-\pi}^\pi h(\theta)d\theta$.
For example, if $f(z)=\sin\sigma z$, then $f\in E_\sigma$
and $D=2\sigma/\pi\approx 0.6366\sigma$
The existence of the limit follows from a theorem of
Levinson \cite{K,L}. Much more precise information
about $n(r)$ under the condition (\ref{2}) is contained
in the theorem of Beurling and Malliavin \cite{BM}.

In the general case, the density might not exist
as was shown by examples
in \cite{KR,R}. Moreover, it is possible that
$D>2\sigma/\pi$, see \cite{EN}.
An easy estimate using Jensen's formula gives
$D\leq 2e\sigma/\pi\approx 1.7305\sigma$. This estimate is exact in
the larger class of entire functions satisfying the condition
$h(\theta)\leq\sigma$, but it is not exact in $E_\sigma$.

In this paper we find the best possible
upper estimate for the upper density of zeros of functions
in $E_\sigma$.
\vspace{.1in}

\noindent
{\bf Theorem.} {\em The upper density
of zeros of a function
$f\in E_\sigma$ does not exceed $c\sigma$
where $c\approx 1.508879$
is the unique solution of the equation
\begin{equation}
\label{eq}
\log({\sqrt{c^2+1}+c})
=\sqrt{1+c^{-2}},\quad\mbox{on}
\quad (0,+\infty).\end{equation}
For every $\sigma>0$
 there exist entire functions $f\in E_\sigma$
such that $D=c\sigma$.
}
\vspace{.1in}

{\em Proof.} Without loss of generality
we assume that $\sigma=1$. Moreover,
it is enough to consider only even functions.
To make a function $f$ even we replace it by
$f(z)f(-z)$, which results in multiplication of both
the indicator $h$ and the upper density $D$
by the same factor
of $2$.

Let $t_n\to+\infty$ be such sequence
that $\lim n(t_n)/t_n=D.$
Consider the sequence of subharmonic functions
$v_n(z)=t_n^{-1}\log|f (t_nz)|$.
Compactness Principle for subharmonic functions
\cite[Theorem 4.1.9]{Hor}
implies that one can choose
a subsequence that converges in $\D'$
(Schwartz's distributions). The limit function $v$
is subharmonic in the plane,
and satisfies
\begin{equation}
\label{u}
v(z)\leq |\Ima z|,\quad z\in\C,\quad\mbox{and}\quad v(0)=0.
\end{equation}
Let $\mu$ be the Riesz measure of this function.
We have to show that
\begin{equation}
\label{haveto}
\mu(\{ z:|z|\leq 1\})\leq c.
\end{equation}
First we reduce the problem to the case that
the Riesz measure $\mu$ is supported by the real line.
We have
$$v(z)=
\frac{1}{2}\int\log\left|1-\frac{z^2}{\zeta^2}\right|d\mu_\zeta.$$
Let us compare this with
$$v^*(z)=\frac{1}{2}\int_0^\infty\log\left|1-\frac{z^2}{t^2}\right|
d\mu^*_t,$$
where $\mu^*$ is the radial projection
of the measure $\mu$: it is supported on $[0,+\infty)$ and
$\mu^*(a,b)=\mu(\{ z:a<|z|<b\}),\; 0\leq a<b.$
It is easy to see that
\begin{equation}
\label{star}
v^*(z)\leq \sigma'|\Ima z|,\quad z\in\C,\quad\mbox{and}\quad
v^*(0)=0
\end{equation}
with some $\sigma'>0.$
We claim that one can choose $\sigma'\leq 1$ in (\ref{star}).
Let $\sigma'$ be the smallest number for which (\ref{star})
holds.
Then, by the subharmonic version of the
theorem of Levinson mentioned above (see, for example,
\cite{MS}),
the limit
$$\lim_{r\to\infty}r^{-1}v^*(rz)=\sigma'|\Ima z|$$
exists in $\D'$ and thus
$$\lim_{r\to\infty}\frac{1}{r}\int_0^r\frac{n_{v^*}(t)}{t}dt
=\lim_{r\to\infty}\frac{1}{2\pi r}\int_{-\pi}^\pi
v^*(re^{i\theta})d\theta=2\sigma'/\pi,$$
where
\begin{equation}
\label{n}
n_{v^*}(r)=\mu^*[0,r]=\mu\{ z:|z|\leq r\}.
\end{equation}
Similar limits exist for $v$, and we have $n_v=n_{v^*}$,
from which we conclude that $\sigma'\leq 1$.

From now on we assume that $v$
is harmonic in the upper and lower half-planes, and that
\begin{equation}\label{y}
v(iy)\sim y,\quad y\to+\infty.
\end{equation}

Let $u$ be the harmonic function in the upper half-plane
such that
$\phi=u+iv$ is analytic, and $\phi(0)=0$.
Then $\phi$ is a conformal map
of the upper half-plane onto some region $G$
of the form
\begin{equation}\label{G}
G=\{ x+iy: y>g(x)\},
\end{equation}
where $g$ is an even upper semi-continuous function,
$g(0)=0$. Moreover,
\begin{equation}\label{as}
\phi(iy)\sim iy,\quad\mbox{as}\quad
y\to+\infty,
\end{equation}
which follows from (\ref{y}),
and
\begin{equation}\label{sim}
\phi(-\overline{z})=-\overline{\phi(z)},
\end{equation}
because both the region $G$ and the normalization of $\phi$
are symmetric with respect to the imaginary axis.
Finally we have
\begin{equation}
\mu([0,x])=\frac{2}{\pi}u(x).
\end{equation}
For all these facts we refer to \cite{Lev}.
\vspace{.1in}

{\em Remark.} The function $\Rea\phi(x)=u(x)$ might be discontinuous for $x\in\R$. We agree to understand
$u(x)$ as the limit from the right $u(x+0)$ which always
exists since $u$ is increasing.
\vspace{.1in}

Inequality (\ref{u}) implies that
$v(x)\leq 0$, thus $g(x)\leq 0$, in other words, $G$
contains the upper half-plane.
\vspace{.1in}

Thus we obtain the following extremal problem: {\em
Among all univalent analytic functions $\phi$
satisfying (\ref{sim})
and mapping the upper half-plane
onto regions of the form (\ref{G})
with $g\leq 0, \; g(0)=0$
and satisfying $\phi(0)=0$ and (\ref{as}),
maximize $\Rea \phi(1).$}
\vspace{.1in}

We claim that the extremal function $g$ for this problem is
$$g_0(x)=\left\{\begin{array}{ll}-\infty,&0<|x|<\pi c/2,\\
                                  0,&\mbox{otherwise},
\end{array}\right.$$
where $c>1$ is the solution of equation (\ref{eq}).
The corresponding region is
shown in Fig. 1.
For the extremal function we have $\phi_0(1)=
\pi c/2-i\infty$.

\vspace{.1in}
\begin{center}
\epsfxsize=3.0in
\centerline{\epsffile{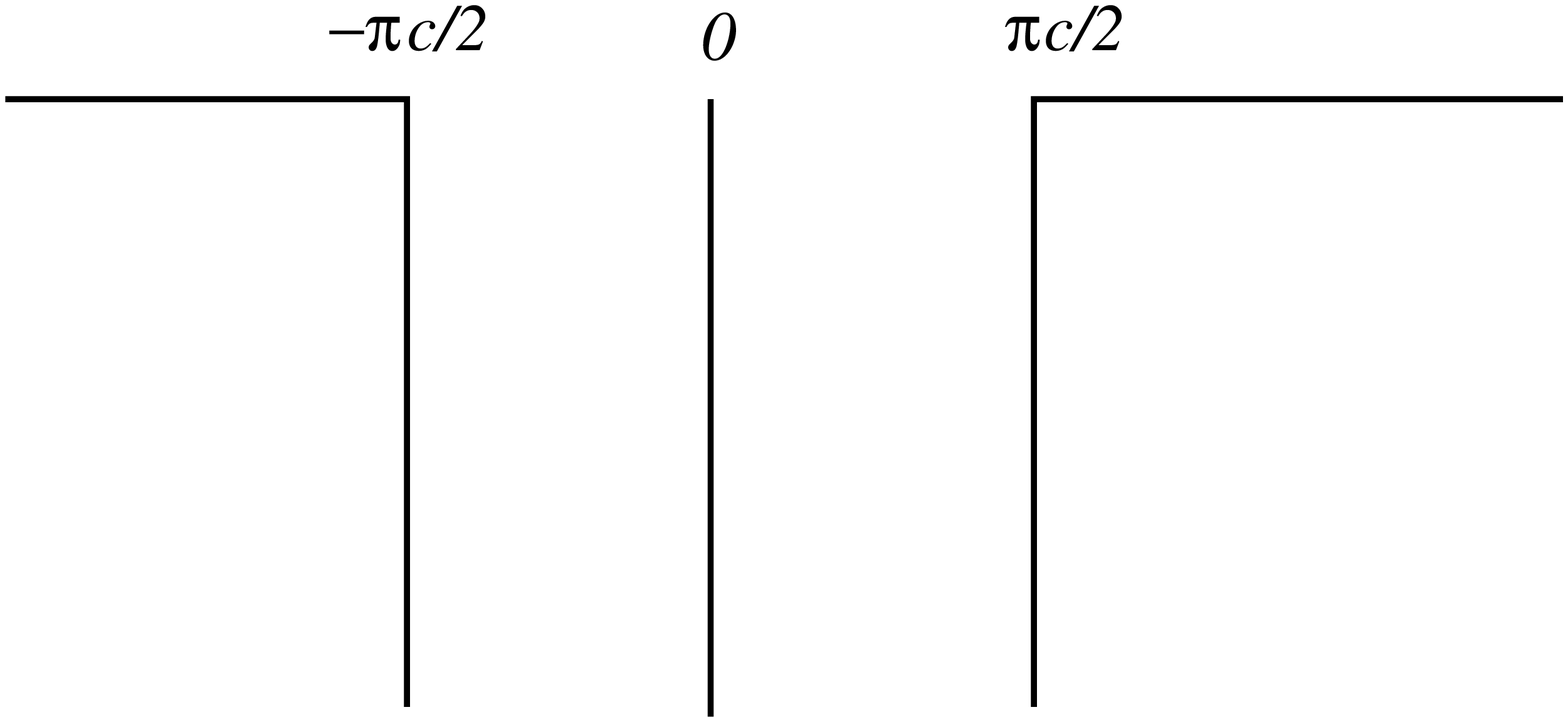}}
\nopagebreak
\vspace{.1in}
Fig. 1. Extremal region.
\vspace{.1in}
\end{center}

To prove the claim,
we fist notice that for a given $G$ the mapping function
is uniquely defined. Let $a=\phi(1)$, and $b=\Rea a$.
Next we show that making $g$ smaller on the interval
$(0,b)$ results in increasing $\Rea\phi(1)$
and making $g$ larger on the interval
$(b,+\infty)$ also results in increasing $\Rea \phi(1)$.
The proofs of both statements
are similar. Suppose that $g_1\leq g,\;
g_1\neq g,$ and $g_1(x)=g(x)$ outside
of the two intervals $p<|x|<q$, where $0<p<q<b$.
Let $G_1$ be the region above the graph of $g_1$,
and $\phi_1$ the corresponding mapping
function normalized in the same
way as $g$. Then $G\subset G_1$,
and the conformal map $\phi_1^{-1}\circ \phi$
is defined in the upper half-plane and maps it into itself.
We have
$$\phi_1^{-1}\circ\phi(x)=x+2x\int_0^\infty
\frac{w(t)}{t^2-x^2}dt,$$
where $w\neq 0$ is a non-negative function supported
on some interval inside $(0,1)$.
Putting $x=1$ we obtain
$$\phi^{-1}_1(a)=1+2\int_0^\infty\frac{w(t)}{t^2-1}dt,$$
so $\phi^{-1}_1(a)<1$, that is $\Rea\phi_1(1)>b$.
This proves our claim.

It remains to compute the constant $b$ in the extremal
domain. We recall that $\phi_0(1)=b-i\infty$ and
assume that $b=\phi_0(k)$ for some $k>1$. Here $\phi_0$
is the extremal mapping function.
Then by the Schwarz--Christoffel formula we have
\begin{equation}\label{int}
\phi_0(z)=\frac{1}{2}
\int_0^{z^2}
\frac{\sqrt{\zeta-k^2}}{\zeta-1}d\zeta.
\end{equation}
To find $k$, we use the condition that
$$\Ima p. v. \int_0^{k^2}
\frac{\sqrt{\zeta-k^2}}{\zeta-1}d\zeta=0.$$
Denoting $c=\sqrt{k^2-1}$ and evaluating the integral,
we obtain
$$\log({\sqrt{c^2+1}+c})=\sqrt{1+c^{-2}}.$$
Finally the jump of the real part of the integral
in (\ref{int}) occurs at the point $1$ and has magnitude
$\pi\sqrt{k^2-1}=\pi c$. This completes the
proof of the upper estimate in Theorem~1.

To construct an example showing that this
estimate can be attained, we follow the construction in \cite[Sect.9-10]{EN}.
The role of the subharmonic function $u_1$ there is played now by our extremal
function $v_0=\Ima\phi_0$.

\vspace{.2in}

{\em Purdue University

West Lafayette IN 47907 USA}

\bigskip

{\em J. Kepler University

Linz A 4040 Austria}

\end{document}